\documentclass[12pt]{amsart}
\usepackage{amsmath,amsfonts}
\usepackage{amssymb}

\newtheorem{theorem}{Theorem}[section]
\newtheorem{proposition}[theorem]{Proposition}
\newtheorem{lemma}[theorem]{Lemma}

\begin{document}

\title[Poisson integrals on non-tube  domains]{Characterization of
  Poisson integrals for non-tube bounded symmetric domains}

\author{Abdelhamid Boussejra}
\address{Abdelhamid Boussejra -- D{\'e}partement de  Math{\'e}matiques, Facult{\'e}
  des sciences, Universit{\'e} Ibn Tofail, K{\'e}nitra, Maroc}
\email{Boussejra@mailcity.com}

\author{Khalid Koufany}
\address{Khalid Koufany --
Institut {\'E}lie Cartan, 
UMR 7502, Universit{\'e} Henri Poincar{\'e}
    (Nancy 1) B.P. 239, 
F-54506 Vand{\oe}uvre-l{\`e}s-Nancy,
France
}
\email{khalid.koufany@iecn.u-nancy.fr}

\begin{abstract} We characterize the $L^p-$range, $1<p<+\infty$, of the
  Poisson transform
  on the Shilov boundary for non-tube bounded symmetric domains. We prove that this
  range is a Hua-Hardy type space for harmonic functions satisfying a Hua system.
\end{abstract}


\maketitle

\section{Introduction}
Let $\Omega=G/K$ be a Riemannian symmetric space of non-compact type. To each boundary $G/P$ one can define a
Poisson transform, which is an integral operator from hyperfunctions on
$G/P$ into the space of eigenfunctions on $\Omega$ of the algebra
$\mathcal{D}(\Omega)^G$ of invariant differential operators.
 For the maximal
boundary,  $G/P_{\min}$,  the most important result is the Helgason
conjecture, proved by Kashiwara {\it et al.} \cite{K et Al.} which states
that a function is eigenfunction of all invariant differential operators on $\Omega$   if and only if it is  Poisson integral 
$$\mathcal{P}_\lambda f(gK)=\int_K f(k)e^{-\langle \lambda+\rho,
  H(g^{-1}k\rangle } dk.$$
  of a
hyperfunction on the maximal boundary, for a generic $\lambda\in\mathfrak{a}_\mathbb{C}^*$.\\
 For other function spaces such as $L^p(G/P_{\min})$ the characterization
 is in connection with Fatou's theorems. We mention here the work of
 Helgason \cite{Helgason1} and Michelson \cite{Michelson} for
 $p=\infty$, and Sj{\"o}rgen \cite{Sjorgen} for $1\leq p<\infty$ using weak
 $L^p-$spaces. Another characterization for $1\leq p\leq \infty$, using Hardy-type spaces, was done
 by Stoll \cite{Stoll} in the harmonic case and by  Ben Sa{\"\i}d {\it et al.}
 \cite{Bensaid_al} in the general case.\\

If $\Omega$ is a bounded symmetric domain,
one is interested in functions whose boundary values are supported
on the Shilov boundary (minimal boundary) $S:=G/P_{\max}$ rather that the maximal boundary $G/P_{\min}$. 
For the Shilov boundary the Poisson transform  is defined by
$$\mathcal{P}_s f(gK)=\int_K f(k)e^{-\langle s\rho_0+\rho_1,
  H_1(g^{-1}k\rangle } dk, \;\;  s\in\mathbb{C}.$$ 
  
In this case, Hua \cite{Hua}
had proved 
 that the algebra of invariant differential operators is not
necessarily the most appropriate for characterizing harmonic functions (i.e., annihilated by
the algebra $\mathcal{D}(\Omega)^G$). Johnson and
Kor{\`a}nyi \cite{Johnson-Koranyi}, generalizing the earlier work of Hua, Kor{\'a}nyi and Stein
\cite{Koranyi-Stein}, and Kor{\'a}nyi and Malliavin \cite{Koranyi-Malliavin}, introduced an invariant
second order ($\mathfrak{k}_{\mathbb{C}}-$valued) operator $\mathcal{H}$, called since, second-order Hua operator
(or Hua system). They showed, in the tube case, that  a function is annihilated by
the Hua operator
if and only if it is the Poisson integral  $\mathcal{P}_{s_0} f$ ($s_0=n/r$)
 of a hyperfunction on the Shilov
boundary. Thus, in the tube case, The Hua operator plays the same role with
respect to the Shilov boundary as the algebra $\mathcal{D}(\Omega)^G$ does
with respect to the maximal boundary. In his paper
\cite{Lassalle},  Lassalle showed the existence of a smaller system (a projection of
the Hua operator) with the same
properties.\\ 

Later Shimeno \cite{Shimeno} generalize the result of Johnson
and Kor{\'a}nyi; namely he proved that a function is  eigenfunction of
$\mathcal{H}$  if and only if it is a Poisson transform $\mathcal{P}_sf$ of
a hyperfunction on the Shilov boundary for generic $s\in\mathbb{C}$. \\

In \cite{Boussejra}, the first author gave a characterization of the
Poisson transform $\mathcal{P}_s$ on $L^p(S)$, which  closes the tube type symmetric domains case characterization.\\

It thus arises the question of characterizing the range of the Poisson
transform $\mathcal{P}_s$
on $L^p(S)$, $1<p<+\infty$, for  non-tube bounded symmetric domains on $L^p(S)$. 
The purpose of this paper is to answer this question.\\

For general bounded symmetric domains the Poisson integrals are not eigenfunctions of
the second-order Hua operator $\mathcal{H}$, see for instance \cite{Berline-Vergne} or
\cite{Koufany-Zhang}. However for type $\mathbf{I}_{r,r+b}$ domains of non-tube type, (see
\cite{Berline-Vergne} and \cite{Koufany-Zhang}) there is a variant of the
second-order Hua operator, $\mathcal{H}^{(1)}$, by
taking the first component of $\mathcal{H}$, since in this case $\mathfrak{k}_\mathbb{C}$
 is a sum of two irreducible ideals $\mathfrak{k}_\mathbb{C}=\mathfrak{k}_\mathbb{C}^{(1)}\oplus
\mathfrak{k}_\mathbb{C}^{(2)}$. It is proved, in \cite{Koufany-Zhang} (and
in \cite{Berline-Vergne} for the harmonic case, $s=(2r+b)/r$) that a smooth function $f$
on $\mathbf{I}_{r,r+b}$ is a solution of the Hua system,
$\mathcal{H}^{(1)}f=\frac{1}{4}(s^2-(r+b)^2)fI_r$ if and only if it is the
Poisson transform $\mathcal{P}_s$ of a hyperfunction on the Shilov boundary.\\

For general non-tube domains, and for the harmonic case (
i.e., for $s=n/r$ in our parametrization)  the characterization of the
image of the Poisson transform $\mathcal{P}_{\frac{n}{r}}$ on
hyperfunctions over the Shilov boundary  was done by Berline and Vergne
\cite{Berline-Vergne}  where certain third-order Hua operator was
introduced. Recently, the second author and Zhang \cite{Koufany-Zhang} generalize the result of
of Berline and Vergne to any (generic) $s$.  They introduce two third-order
Hua operators $\mathcal{U}$ and $\mathcal{W}$ (different from the Berline
and Vergne operator) and prove that an eigenfunction $f$ of
$\mathcal{D}(\Omega)^G$ is a solution the Hua system (\ref{eq:3-rd-cha})
if and only if it is a Poisson transform of a hyperfunction on the Shilov
boundary.\\

Let $\mathcal{E}_s(\Omega)$ be the space of harmonic functions on $\Omega$
that are solutions of the Hua system (for type $\mathbf{I}_{r,r+b}$
domains, an eigenfunction of $\mathcal{H}^{(1)}$ is indeed harmonic). Then
the image $\mathcal{P}_s(L^p(S))$  is a proper closed
subspace of $\mathcal{E}_s(\Omega)$. For $1<p<+\infty$, we introduce the Hua-Hardy type space,
$\mathcal{E}^p_s(\Omega)$ of functions $f\in \mathcal{E}_s(\Omega)$ such
that
$$\|f\|_{s,p}=\sup_{t>0}e^{-t(\Re(s)r-n)}(\int_K
|f(ka_t)|^pdk)^{1/p}<+\infty.$$
Our main result (see Theorem \ref{second_main_Thm}) says that {\it if
$s\in\mathbb{C}$ is such that $\Re(s)>\frac{a}{2}(r-1)$, a smooth function
$F$ on $\Omega$ is the Poisson transform $F=\mathcal{P}_sf$ of a function
$f\in L^p(S)$ if and only if $f\in\mathcal{E}^p_s(\Omega)$}. Our method of
proving this characterization uses  an $L^2$ version of this theorem (see Theorem \ref{main_th_L_2}) and an
inversion formula for the Poisson transform (see Proposition
\ref{inversion_L_2}) which needs Fatou-type
theorems (see Theorem \ref{thm_Fatou1} and Theorem \ref{thm4.5}).\\





\section{Preliminaries}\label{prelim}
Let $\Omega$ be an irreducible bounded symmetric domain in a complex
$n-$dimensional space $V$. Let $G$ be be the identity
component of the group of biholomorphic
automorphisms of $\Omega$, and $K$ be the isotropy subgroup of $G$ at the
point $0\in \Omega$. Then $K$ is a maximal compact subgroup of $G$ and as a
Hermitian symmetric space, $\Omega=G/K$. 
Let $\mathfrak g$ be the Lie algebra of $G$, and 
$$\mathfrak{g}=\mathfrak{k}+\mathfrak{p}$$
 be its Cartan decomposition.
The Lie algebra $\mathfrak{k}$ of $K$ has one dimensional center 
$\mathfrak{z}$. Then there
exists an element $Z_0\in\mathfrak{z}$ such that $\textrm{ad} {Z_0}$ defines
the complex structure of $\mathfrak{p}$. 
Let
\begin{equation*}
\mathfrak{g}_\mathbb{C}=\mathfrak{p}^+\oplus\mathfrak{k}_\mathbb{C}\oplus\mathfrak{p}^-
\end{equation*}
be the corresponding eigenspaces decomposition of 
$\mathfrak{g}_\mathbb{C}$, the complexification of $\mathfrak{g}$. Let
$G_\mathbb{C}$ be a connected Lie group with Lie algebra
$\mathfrak{g}_\mathbb{C}$ and $P^+$, $K_\mathbb{C}$, $P^-$ be the analytic
subgroups of $G_\mathbb{C}$ corresponding to  $\mathfrak{p}^+$,
$\mathfrak{k}_\mathbb{C}$, $\mathfrak{p}^-$. Denote by $\sigma$ the
conjugaison of $G_\mathbb{C}$ with respect to $G$. Then we have,
$\sigma(P^{\pm})\subset P^\mp$ and $\sigma(K_\mathbb{C})\subset
K_\mathbb{C}$.\\

Let $\mathfrak{h}$ be a maximal Abelian subalgebra of $\mathfrak{k}$, and
let $\Delta(\mathfrak{g}_\mathbb{C},\mathfrak{h}_\mathbb{C})$ be the
corresponding set of roots. As $Z_0$ belongs to $\mathfrak{h}$, the space
$\mathfrak{p}^+$ is stable by $\text{ad}\mathfrak{h}$. The roots
$\gamma\in\Delta(\mathfrak{g}_\mathbb{C},\mathfrak{h}_\mathbb{C})$ such that
$\mathfrak{g}^\gamma\subset \mathfrak{p}^+$ are said to be positive
non-compact, and we denote by $\Phi$ the set of such roots. Let $\gamma\in
\Phi$, then one may choose elements $H_\gamma\in i\mathfrak{h}$, $E_\gamma\in
\mathfrak{g}^\gamma$, $E_{-\gamma}\in
\mathfrak{g}^{-\gamma}$ such that $[E_\gamma,E_{-\gamma}]=H_\gamma$ and
$\sigma(E_\gamma)=-E_{-\gamma}$.
 Let $X_\gamma=E_\gamma+E_{-\gamma}$ and
 $Y_\gamma=i(E_\gamma-E_{-\gamma})$. Then, by a classical Harish-Chandra
 construction, there exists a maximal set $\Gamma=\{\gamma_1,\ldots, \gamma_r\}$
  of strongly orthogonal roots in  $\Phi$. For simplicity, let us set for,
  $1\leq j\leq r$, 
$$E_j=E_{\gamma_j},\;\;  X_j=X_{\gamma_j},\;\; Y_j=Y_{\gamma_j}.$$
Then,
$$\mathfrak{a}=\sum_{j=1}^r\mathbb{R}X_j,$$
is a Cartan subspace of the pair $(\mathfrak{g},\mathfrak{k})$. Let
$\mathfrak{a}^*$ denote the dual of $\mathfrak{a}$ and let $\{\beta_1,
\beta_2,\ldots, \beta_r\}$ be a basis of $\mathfrak{a}^*$ determined by
$$\beta_j(X_k)=2\delta_{j,k},\quad 1\leq j, k\leq n.$$
 The restricted root system
$\Sigma=\Sigma(\mathfrak{g},\mathfrak{a})$ of $\mathfrak{g}$
relative to $\mathfrak{a}$ is (of type $C_r$ or $BC_r$) given by
$$ \pm \beta_j\;\;\; (1\leq j\leq r)\quad \text{each with multiplicity }\;
1,$$
$$ \pm \frac{1}{2}(\beta_j\pm\beta_k)\;\;\; (1\leq j\not=k\leq r)\quad
\text{each with multiplicity }\; a,$$
and possibly
$$ \pm \frac{1}{2}\beta_j\;\;\; (1\leq j\leq r)\quad \text{each with multiplicity }\;
2b.$$
Let $\Sigma^+=\{\beta_j,\frac{1}{2}\beta_j,
\frac{1}{2}(\beta_\ell\pm\beta_k) ;\;\; 1\leq j\leq r,\; 1\leq \ell\not=
k\leq r\}$ the set of positive restricted roots.
Then the set
$\Lambda=\{\alpha_1,\; \ldots,\; \alpha_{r-1},\;
\alpha_r\}$ of
simple roots in $\Sigma^+$ is such that 
\begin{equation*}
\alpha_j=\frac{1}{2}(\beta_{r-j+1}-\beta_{r-j}),\; 1\leq j\leq r-1
\end{equation*}and
\begin{equation*}
\alpha_r=\begin{cases} \beta_1 & \text{for tube case}\\ \frac{1}{2}\beta_1 & \text{for non-tube case}.\end{cases}
\end{equation*}
Let $\Lambda_1=\{\alpha_1,\ldots,\alpha_{r-1}\}$ and write
$\Sigma_1=\Sigma\cap\mathbb{Z}\cdot \Lambda_1$. Define
$$\mathfrak{m}_{1,1}=
\mathfrak{m} +\mathfrak{a} +\sum_{\gamma\in
  \Sigma_1}\mathfrak{g}^\gamma,\;\; 
\mathfrak{n}_1^+=\sum_{\gamma\in
  \Sigma^+\setminus\Sigma_1}\mathfrak{g}^\gamma.$$
Let 
$$\mathfrak{a}_1=\{H\in \mathfrak{a} \; :\; \gamma(H)=0\; \forall \gamma\in
\Lambda_1 \},$$ 
 then $\mathfrak{m}_{1,1}$ is the
centralizer of $\mathfrak{a}_1$ in $\mathfrak{g}$ and
$\mathfrak{p}_1=\mathfrak{m}_{1,1}+\mathfrak{n}_1^+$ 
is a standard
parabolic subalgebra of $\mathfrak{g}$ with Langlands decomposition
$\mathfrak{m}_1+\mathfrak{a}_1+\mathfrak{n}_1^+$, where $\mathfrak{m}_1$ is
the orthocomplement of $\mathfrak{a}_1$ in $\mathfrak{m}_{1,1}$ with respect
to the Killing form. Note that
$\theta(\mathfrak{n}_1^+)=\sum_{\gamma\in\Sigma^+\setminus\Sigma_1}\mathfrak{g}^{-\gamma}$.
 Let $P_1$
be the corresponding parabolic subgroup and $P_1=M_1 A_1 N^+_1$ its
Langlands decomposition. Obviously, $P_1$ is a maximal parabolic subgroup
of $G$, thus 
the Shilov boundary $S$ can be viewed
as $S=G/P_1=K/K_1$, where $K_1=M_1\cap K$ \\

If we define the element $X_0=\sum_{j=1}^rX_j$, Then
$\mathfrak{a}_1=\mathbb{R}X_0$. Let $$\mathfrak{a}(1)=\sum_{j=1}^{r-1}\mathbb{R}(X_j-X_{j+1})$$ be the
orthocomplement of $\mathfrak{a}_1$ in $\mathfrak{a}$ with respect to the
Killing form,

\begin{equation}\label{orthoproj}
\mathfrak{a}=\mathfrak{a}_1\oplus^\bot \mathfrak{a}(1)
=\mathbb{R}X_0\oplus^\bot\sum_{j=1}^{r-1}\mathbb{R}(X_j-X_{j+1}).
\end{equation}

We denote $\rho_0$ the linear form on $\mathfrak{a}_1$ such that, $\rho_0(X_0)=r$. We extend $\rho_0$ to
$\mathfrak{a}$ via the orthogonal projection (\ref{orthoproj}).
If $\rho_1$ is the restriction of $\rho$ to $\mathfrak{a}_1$, then it is clear
that
$$\rho_1(X_0)=rb+r+a\frac{r(r-1)}{2}=n.$$
Again, we extend $\rho_1$ to $\mathfrak{a}$ via the orthogonal
projection (\ref{orthoproj}). Then
$$\rho_1=(b+1+a\frac{(r-1)}{2})\rho_0=\frac{n}{r}\rho_0.$$
For $g\in G$, define $H(g)\in\mathfrak{a}$ as the unique element such that
\begin{equation*}
g\in K\exp(H(g))N\subset KAN=G.
\end{equation*} We also denote by $\kappa(g)\in K$ and  $H_1(g)\in \mathfrak{a}_1$
the unique elements such that 
\begin{equation*}
g\in \kappa(g)M_1\exp(H_1(g))N_1\subset KM_1A_1N_1=G.
\end{equation*}

The following lemma will be useful for the sequel. 
\begin{lemma}[{\cite[Lemma 6.1.6]{Schli}}]\label{lemme_proj}
\begin{itemize}
\item[$(i)$] Let $x, y\in G$, $\bar{n}\in\bar{N}_1$ and $a\in A_1$. Then
\begin{eqnarray}
H_1(x\kappa(y))=H_1(xy)-H_1(y)\label{lemme_proj1}\\
H_1(\bar{n}a^{-1})=H_1(\bar{n})-H_1(a)\label{lemme_proj2}
\end{eqnarray}
\item[$(ii)$] Let $t>0$ and $\bar{n}\in \bar{N}_1$. Then 
\begin{equation}\label{lemme_proj3}
  \rho_0(H_1(a_t\bar{n}a_{-t}))\leq \rho_0(H_1(\bar{n})).
\end{equation}
\end{itemize}
\end{lemma}


\section{The Poisson transform and the Hua operators}\label{sect_Hua_Poisson}
For any real analytic manifold $X$, we denote by $\mathcal{B}(X)$ the space
of all hyperfunctions on $X$. 
We will view a function on the Shilov boundary $S=G/P_1$ as a $P_1-$invariant
function on $G$. For $s\in\mathbb{C}$, we denote by
$\mathcal{B}(G/P_1;s)$ the space of hyperfunctions $f$ on $G$ satisfying
\begin{equation*}
f(gman)=e^{(s\rho_0-\rho_1)\log a}f(g),\;\;\forall g\in G,\; m\in
M_1,\; a\in A_1,\; n\in N_1^+,
\end{equation*}
 The Poisson transform 
of a function $f\in \mathcal{B}(G/P_1;s)$, is defined by
\begin{equation*}
\mathcal{P}_s f(gK)=\int_K e^{-\langle s\rho_0+\rho_1,H_1(g^{-1}k) \rangle}f(k)dk.
\end{equation*}
Since $G=KP_1$, the restriction from $G$ to $K$
defines a $G-$isomorphism from $\mathcal{B}(G/P_1, s)$ onto the space
$\mathcal{B}(K/K_1)$ of all hyperfunctions $f$ on $K$ such that
$f(kh)=f(k)$ for all $h\in K_1$.\\

We review the construction of Hua operators of the second order (see \cite{Johnson-Koranyi}) and the third
order (see \cite{Berline-Vergne}, \cite{Koufany-Zhang}).  \\
Let $\{v_j\}$ be 
a basis of $\mathfrak{p}^+$ and $\{v^*_j\}$ be th dual basis of
$\mathfrak{p}^-$ with respect to the Killing form. 
Let
$\mathcal{U}(\mathfrak{g}_\mathbb{C})$ denote the universal enveloping
algebra of $\mathfrak{g}_\mathbb{C}$. The second-order Hua operator, is the
element of $\mathcal{U}(\mathfrak{g}_\mathbb{C})\otimes
\mathfrak{k}_\mathbb{C}$ defined by
\begin{equation*}
\mathcal{H}=\sum_{i,j}v_iv^*_j\otimes[v_j,v^*_i]
\end{equation*}
It is known that the Hua operator does not depends on basis, therefor, for
computations one can choose the root vectors basis $\{E_j\}_{j=1}^r$.\\


For tube
domains the Hua operator $\mathcal{H}$ maps the Poisson kernels $$P_s(gk)=e^{-\langle s\rho_0+\rho_1,H_1(g^{-1})\rangle}$$
into the center of 
$\mathfrak k_{\mathbb C}$, namely the Poisson kernels
are  its eigenfunctions up to an  element in the center,
but it is not true for non-tube domains, see \cite[Theorem
5.3]{Koufany-Zhang}.
However for non-tube type $\mathbf{I}$ domains, $\mathbf{I}_{r,r+b}\simeq
SU(r,r+b)/S(U(r)\times U(r+b))$, the situation is not quite different from
the tube case. 
In fact, $\mathfrak k_{\mathbb C}$ is a sum of two irreducible ideals,
$\mathfrak k_{\mathbb C}=\mathfrak k_{\mathbb C}^{(1)}+\mathfrak
k_{\mathbb C}^{(2)}$ where
$$\begin{array}{ll}
\mathfrak k_{\mathbb C}^{(1)}&=\left\{\begin{pmatrix}A&0\\
    0&\frac{\textrm{tr}(A)}{r+b}I_{r+b}\end{pmatrix},\; A\in\mathfrak{gl}(r+b,\mathbb{C})\right\},\\
\mathfrak k_{\mathbb
  C}^{(1)}&=\left\{\begin{pmatrix}0&0\\0&D\end{pmatrix},\; D\in\mathfrak{sl}(r+b,\mathbb{C})\right\}.
\end{array}$$
There is a variant of the Hua
operator, $\mathcal{H}^{(1)}$ see \cite{Berline-Vergne}, \cite{Koufany-Zhang},  by taking  the 
projection of $\mathcal{H}$ onto $\mathfrak{k}_{\mathbb{C}}^{(1)}$. In
\cite{Koufany-Zhang}, the second author and Zhang  
showed that the operator $\mathcal{H}^{(1)}$ has the Poisson kernels
as its eigenfunctions and they  found
the eigenvalues.
They proved further  that the eigenfunctions of the Hua operator
$\mathcal{H}^{(1)}$ are harmonic functions (i.e., eigenfunctions of all invariant differential operators on $\Omega$),
and gave the following characterization
of the range of the Poisson transform 
for  $\mathbf{I}_{r,r+b}$.
\begin{theorem}[{\cite[Theorem 6.1]{Koufany-Zhang}}]\label{KZ1}
Suppose $s\in\mathbb{C}$ satisfies the 
following condition 
\begin{equation*}
-4[b+1 +j +\frac{1}{2}(s-r-b)]\notin \{1, 2, 3, \cdots\},\; \; \text{for\/}\;\; j=0
\;\text{and\/}\; 1.
\end{equation*}
Then the Poisson transform $\mathcal{P}_s$ is a $G-$isomorphism of
$\mathcal{B}(S)$ onto the space of smooth functions
$f$ on
$\Omega$ that satisfy
\begin{equation} \label{eq:hua_typeI}
\mathcal{H}^{(1)}f=\frac{1}{4}(s^2-(r+b)^2)fI_r. 
\end{equation}
\end{theorem}

For the characterization
of range of the Poisson transform for general 
non-tube domains the second author and Zhang \cite{Koufany-Zhang} introduced new
third-order Hua operators $\mathcal{U}$ and $\mathcal{W}$ :
\begin{equation*}
\mathcal{U}=\sum_{i,j,k}v^*_iv^*_jv_k\otimes[v_i,[v_j,v^*_k]],
\end{equation*}

\begin{equation*}
\mathcal{W}=\sum_{i,j,k}v_kv^*_iv_j\otimes[[v^*_k,v_i],v_j]],.
\end{equation*}
Similarly to $\mathcal{H}$, the operators $\mathcal{U}$ and $\mathcal{W}$
do not depend on the basis. \\

Denote
$$c=2(n+1)+\frac{1}{n}(a^2-4)\dim(\mathcal{P}^{(1,1)}),$$
where $\mathcal{P}^{(1,1)}$ is the dimension of the irreducible subspaces
of holomorphic polynomials on $\mathfrak{p}^+$ with lowest weight
$-\gamma_1-\gamma_2$. For any $s\in\mathcal{C}$, put
$\sigma=\frac{1}{2}(s+\frac{n}{r})$. For general non-tube domains we have
the following

\begin{theorem}[{\cite[Theorem 7.2]{Koufany-Zhang}}]\label{KZ2}
Let $\Omega$
be a  bounded symmetric non-tube domain of rank $r$ in
$\mathbb{C}^n$. Suppose $s\in\mathbb{C}$ satisfies 
\begin{equation*}
-4[b+1 +j\frac{a}{2} +\frac{1}{2}(s-\frac{r}{n})]\notin \{1, 2, 3, \cdots\},\; \; \text{for\/} \;\; 
j=0 \; \text{and\/} \;  1.
\end{equation*}
Then the Poisson transform $\mathcal{P}_s$ is a $G-$isomorphism of
$\mathcal{B}(S)$ onto the space of harmonic
functions $f$ on $\Omega$ that satisfy
\begin{equation}
  \label{eq:3-rd-cha}
\left(\mathcal{U}-\frac{-2\sigma^2+2p\sigma+c}{\sigma(2\sigma-p-b)}
\mathcal{W}\right)
f=0,
\end{equation}
\end{theorem}


\section{The $L^p-$range of the Poisson transform}
For $1<p<+\infty$, we will consider the space $L^p(S)=L^p(K/K_1)$ as the space of all complex
valued measurable (classes) functions $f$ on $K$ that are $K_1-$invariant
and satisfying
$$\|f\|_p=\left(\int_K|f(k)|^pdk\right)^{1/p}<+\infty,$$
where $dk$ is the Haar measure of $K$. 
Let $d\bar{n}$ be
the invariant measure on $\bar{N}_1=\theta(N_1)$ with the normalization 
\begin{equation}\label{normalisation_int}
\int_{\bar{N}_1}e^{\langle -2\rho_1,H_1(\bar{n})\rangle}d\bar{n}=1.
\end{equation}
Then for a continuous function $f$ on $S$ we have
\begin{equation}\label{formule_KtoN}
\int_Kf(k)dk=\int_{\bar{N}_1}f(\kappa(\bar{n}))e^{-2\langle \rho_1,H_1(\bar{n})\rangle}d\bar{n}.
\end{equation}
The space $L^p(S)$ can be viewed as a subspace of $\mathcal{B}(S)$,
thus its image $\mathcal{P}_s(L^p(S))$ is a proper closed subspace of
$\mathcal{E}_s(\Omega)$. We will now, for specific $s$, characterize this
image. For this we need some information on the  integral $\mathbf{c}_s$ in
following proposition.

\begin{proposition}\label{convergence_c_s}
For $s\in\mathbb{C}$ such that $\Re(s)>\frac{a}{2}(r-1)$, the integral
\begin{equation*}
\mathbf{c}_s=\int_{\bar{N}_1}e^{-\langle s\rho_0+\rho_1,H_1(\bar{n})\rangle}d\bar{n}
\end{equation*}
converges absolutely to a constant $\mathbf{c}_s\not=0$.
\end{proposition}
\begin{proof}
For $s\in\mathbb{C}$ let $\lambda_s\in\mathfrak{a}_\mathbb{C}^*$ be the
linear form defined by
\begin{equation*}
\lambda_s(H)=(s\rho_0-\rho_1)(H_1) +\rho(H),\;\; H\in\mathfrak{a}
\end{equation*}
where $H_1$ is the projection of $H$ onto $\mathfrak{a}_1$.
Then the condition $\Re(s)>\frac{a}{2}(r-1)$ is equivalent to 
\begin{equation}\label{cd:2}
\Re (\langle \lambda_s,\alpha\rangle) > 0 \;\; \forall \alpha\in\Sigma^+\setminus\Sigma_1.
\end{equation}
Moreover, we can choose (see for example \cite[Lemma 6.1.4]{Schli}) $\omega$ in the Weyl
group $W$ of $\Sigma$ such that
$$\begin{array}{rl}
(i)&\omega\cdot H=H,\;\; \; \forall H\in \mathfrak{a}_1,\\
(ii)& \omega(\Sigma^+\cap\Sigma_1)=-\Sigma^+\cap\Sigma_1,\\
(iii)&\omega(\Sigma^+\setminus\Sigma_1)=\Sigma^+\setminus\Sigma_1.
\end{array}$$
Since $\langle\omega\lambda_s,\alpha \rangle=\langle\lambda_s,\omega^{-1}\alpha \rangle$, the condition (\ref{cd:2}) is equivalent to
\begin{equation*}
\Re(\langle \omega\lambda_s,\alpha\rangle)>0,\;\;\; \forall \alpha\in\Sigma^+
\end{equation*}
Furthermore
$$\langle s\rho_0+\rho_1, H_1(g)\rangle=\langle \omega\lambda_s+\rho, H(g)\rangle$$
so that
$$ \int_{\bar{N}_1}e^{-\langle s\rho_0+\rho_1,H_1(\bar{n})\rangle}d\bar{n}=
\int_{\bar{N}_1}e^{-\langle \omega\lambda_s+\rho,H(\bar{n})\rangle}d\bar{n}$$
and the  right hand side is the Harish-Chandra $c$ function,
$c(\omega\lambda_s)$ associated with the maximal parabolic subgroup, which converges absolutely, see \cite{Schli}.\\
\end{proof}

Let $s\in\mathbb{C}$. Let $\mathcal{E}_s(\Omega)$ be the space of
harmonic functions on $\Omega$ that satisfy (\ref{eq:hua_typeI}) in type
I domains or (\ref{eq:3-rd-cha}) in general domains. 
It is clear that the image $\mathcal{P}_s(L^p(S))$ is a proper closed
subspace of the eigenspace $\mathcal{E}_s(\Omega)$. Hence, it is natural to
look for a characterization of those $F\in\mathcal{E}_s(\Omega)$ that are
Poisson transform of some $f\in L^p(S)$. 

For any $1<p<\infty$, let $\mathcal{E}_s^p(\Omega)$
denote the Hua-Hardy type space of functions $f\in\mathcal{E}_s(\Omega)$ such that
\begin{equation*}
\|f\|_{s,p}=\sup_{a\in A_1}e^{-\langle \Re(s)\rho_0-\rho_1,\log a\rangle }(\int_K|f(ka)|^pdk)^{1/p}<+\infty.
\end{equation*}
Since $\mathfrak{a}_1=\mathbb{R}X_0$, the above integral becomes
\begin{equation*}
\|f\|_{s,p}=\sup_{t>0}e^{-t(\Re(s)r-n)}(\int_K|f(ka_t)|^pdk)^{1/p},
\end{equation*}
where $a_t=\exp(tX_0)$.\\

\subsection{Fatou type theorems}
As a preparation to Fatou-type theorems we prove the following

\begin{proposition}\label{thm_dom}
Let $s\in\mathbb{C}$ be such that $\Re(s)>\frac{a}{2}(r-1)$. Let $\Psi_t$ be
the function defined on $\bar{N}_1$ by
$$\Psi_t(\bar{n})=e^{-\langle s\rho_0+\rho_1, H_1(\bar{n})\rangle+\langle
  s\rho_0-\rho_1,H_1(a_t\bar{n}a_{-t})\rangle}.
$$ Then there exists a non-negative function $\Phi\in L^1(\bar{N}_1)$ such that
$\Psi_t\leq \Phi$  for each $t$.
\end{proposition}

\begin{proof}
It follows from (\ref{lemme_proj3}), that for any $t>0$ and for any $\bar{n}\in\bar{N}_1$,
\begin{equation*}
0\leq \rho_0(H_1(a_t\bar{n}a_{-t})\leq \rho_0(H_1(\bar{n}).
\end{equation*}
Therefore, 
$$|\Psi_t(\bar{n})|\leq 
\begin{cases}
e^{-\langle \Re(s)\rho_0+\rho_1, H_1(\bar{n})\rangle}& \text{if
$\frac{a}{2}(r-1)<\Re(s)\leq \frac{a}{2}(r-1)+b+1$}\\
e^{-2\langle \rho_1, H_1(\bar{n})\rangle}&\text{if $\Re(s)>\frac{a}{2}(r-1)+b+1$}
\end{cases}$$
and the second hand is an integrable function on $\bar{N}_1$ by (\ref{normalisation_int})
and Proposition \ref{convergence_c_s}
.\\
\end{proof}

Let $\mathcal{C}(S)$ be the space of complex-valued continuous functions on
$S$ with the topology of uniform convergence.

\begin{theorem}\label{thm_Fatou1}
Let $s\in\mathbb{C}$ be such that $\Re(s)>\frac{a}{2}(r-1)$. Then
\begin{equation*}
f(k)=\mathbf{c}_s^{-1}\lim_{t\to+\infty}e^{-(rs-n)t}\mathcal{P}_sf(ka_t)
\end{equation*}
 uniformly, for $f\in \mathcal{C}(S)$.
\end{theorem}

\begin{proof}
Let $f\in\mathcal{C}(S)$, then
$$\mathcal{P}_sf(ka_t)=\int_K e^{-\langle s\rho_0+\rho_1,
  H_1(a_{-t}h)\rangle} f(kh)dh.$$
We transform this integral using the formula (\ref{formule_KtoN}) to an integral over $\bar{N}_1$,
\begin{equation*}
\mathcal{P}_sf(ka_t)
=\int_{\bar{N}_1} e^{-\langle s\rho_0+\rho_1,
  H_1(a_{-t}\kappa(\bar{n}))\rangle}f(k\kappa(\bar{n}))e^{-2\langle \rho_1,H_1(\bar{n})
  \rangle}d\bar{n},
\end{equation*}
and by (\ref{lemme_proj1}) we get 
\begin{equation*}
\mathcal{P}_sf(ka_t)
=\int_{\bar{N}_1} e^{-\langle s\rho_0+\rho_1,
  H_1(a_{-t}\bar{n})\rangle} e^{\langle s\rho_0-\rho_1,H_1(\bar{n})
  \rangle}f(k\kappa(\bar{n}))d\bar{n}.\\
\end{equation*}
which by the substitution 
$\bar{n}\mapsto a_{-t}\bar{n}a_t$ and (\ref{lemme_proj2}), becomes 
$$\mathcal{P}_sf(ka_t)=
e^{\langle s\rho_0-\rho_1,H_1(a_t) \rangle} \times \atop
\int_{\bar{N}_1} 
e^{-\langle s\rho_0+\rho_1, H_1(\bar{n})\rangle+\langle
  s\rho_0-\rho_1,H_1(a_t\bar{n}a_{-t})\rangle}
f(k\kappa(a_t\bar{n}a_{-t}))d\bar{n}.$$

But $\rho_1=\frac{n}{r}\rho_0$, and $a_t\bar{n}a_{-t}\to e$ when $t\to
+\infty$, thus, by Proposition \ref{thm_dom},
$$\lim_{t\to+\infty}e^{-(rs-n)t}\mathcal{P}_sf(ka_t)=\mathbf{c}_s f(k).$$
\end{proof}

Let 
$$\varphi_s(a_t):=\int_K e^{-\langle s\rho_0+\rho_1,H(a_{-t}k) \rangle
}dk.$$
then, it follows from the above theorem that
\begin{equation}\label{eq:c_s}
\lim_{t\to\infty}e^{-(rs-n)t}\varphi_s(a_t)=\mathbf{c}_s,\;\;\; \text{if }\; \Re(s)>\frac{a}{2}(r-1).
\end{equation}

As a consequence we can prove the following

\begin{proposition}\label{propo4.4}
Let $s\in\mathbb{C}$ be such that $\Re(s)>\frac{a}{2}(r-1)$. Then there
exists a positive constant $\gamma_s$ such that, for $1<p<\infty$ and $f\in
L^p(S)$, we have
\begin{equation*}
\left(\int_K|\mathcal{P}_s f(ka_t)|^pdk\right)^{1/p}\leq \gamma_s\, e^{(rs-n)t}\|f\|_p.
\end{equation*}
\end{proposition}
\begin{proof}
For $t>0$, we define the function $p^t_s$ by
\begin{equation*}
p^t_s(k)=e^{-\langle s\rho_0+\rho_1,H_1(a_{-t}k^{-1}) \rangle},\;\;\; k\in K.
\end{equation*}
Then the Poisson transform can be written as the convolution
\begin{equation*}
\mathcal{P}_s f(ka_t)=(f\ast p^t_s)(k).
\end{equation*}
Hence, to prove the proposition we use the Haussedorf-Young inequality,
$$\left(\int_K|\mathcal{P}_s f(ka_t)|^pdk\right)^{1/p}\leq \|p^t_s\|_1\,
\|f\|_p \qquad (p>1), $$
and (\ref{eq:c_s}). \\
\end{proof}

Let, as usual, $\hat{K}$, be
the set of equivalence classes of finite dimensional irreducible
representations of $K$.  For $\delta \in \hat{K}$, let
$\mathcal{C}(S)_\delta$ be the linear span of all $K-$finite vectors on $S$
of type $\delta$. It is well known that the space
$\mathcal{C}^K(S):=\oplus_{\delta\in\hat{K}}\mathcal{C}(S)_\delta$ is dense
in $\mathcal{C}(S)$. Recall also, that the space $\mathcal{C}(S)$ is dense
in $L^p(S)$ for $1<p<\infty$.

\begin{theorem}\label{thm4.5}
Let $s\in\mathbb{C}$ be such that $\Re(s)>\frac{a}{2}(r-1)$. Then
\begin{equation*}
f(k)=\mathbf{c}_s^{-1}\lim_{t\to+\infty}e^{-(rs-n)t}\mathcal{P}_sf(ka_t)
\end{equation*}
in $L^p(S)$, for $1<p<\infty$.
\end{theorem}

\begin{proof}
Let $f\in L^p(S)$. By the above density arguments, for any  $\epsilon >0$,  there exists $\varphi\in \mathcal{C}^K(S)$ such that $\|f-\varphi\|_p<\epsilon$. Then we have
$$
\| \mathbf{c}_s^{-1} e^{-(rs-n)t}P^t_sf-f\|_p\leq \| \mathbf{c}_s^{-1} e^{-(rs-n)t}P^t_s(f-\varphi)\|_p
+\atop + \| \mathbf{c}_s^{-1} e^{-(rs-n)t}P^t_s\varphi-\varphi\|_p
+\| \varphi-f\|_p
$$
where the function $P^t_sf$ is defined  by 
\begin{equation}\label{op_P_t}
P^t_sf(k)=\mathcal{P}_sf(ka_t).
\end{equation}  
By Proposition \ref{propo4.4}, 
 
$$\| \mathbf{c}_s^{-1} e^{-(rs-n)t}P^t_s(f-\varphi)\|_p \leq \gamma_s | \mathbf{c}_s^{-1}| \| f-\varphi\|_p,$$
and by Theorem \ref{thm_Fatou1}
$$\lim_{t\to+\infty}\| \mathbf{c}_s^{-1} e^{-(rs-n)t}P^t_s\varphi-\varphi\|_p=0.$$
Thus, $\lim_{t\to+\infty}\| \mathbf{c}_s^{-1} e^{-(rs-n)t}P^t_sf-f\|_p\leq \epsilon(\gamma_s+1)$ and this proves the theorem.
\end{proof}

We can now prove the following estimates
\begin{proposition}\label{nece_cond}
Let $s\in\mathbb{C}$ be such that $\Re(s)>\frac{a}{2}(r-1)$. Then there exists a positive constant $\gamma_s$ such that for $1<p<+\infty$ and $f\in L^p(S)$,
\begin{equation}\label{estimate_L_p}
|\mathbf{c}_s| \|f\|_p\leq \| \mathcal{P}_s f\|_{s,p}\leq \gamma_s\|f\|_p.
\end{equation}
\end{proposition}
\begin{proof} 
In fact, the right hand side of the estimate (\ref{estimate_L_p})  follows from
Proposition \ref{propo4.4}. On the other hand, by Theorem \ref{thm4.5} we have 
$$\lim_{t\to\infty} e^{(n-rs)t}\mathcal{P}_sf(ka_t)=\mathbf{c}_sf(k)$$ in
$L^p(S)$. Hence, there exists a sequence $(t_j)_j$, with $t_j\to+\infty$
when $j\to +\infty$ such that $\lim_{j\to+\infty}e^{(n-rs)t_j}\mathcal{P}_sf(ka_{t_j})=\mathbf{c}_sf(k)$, almost every where in $K$. Therefore, by the classical Fatou lemma, $$|\mathbf{c}_s| \|f\|_p\leq \sup_j e^{(n-r\Re(s))t_j}\| P^{t_j}_s f\|_p$$
and this is how we prove the left hand side of (\ref{estimate_L_p}).
\end{proof}

\subsection{The $\mathrm{L}^2-$Poisson transform range}
Recall that
$$L^2(S)=\oplus_{\delta\in\hat{K}}V_\delta$$where
$V_\delta$ is the finite linear span of
$\{\varphi_\delta\circ k,\;\; k\in K\}$, where $\varphi_\delta$ is the
zonal spherical function corresponding to $\delta$.\\
For $s\in\mathbb{C}$ and $\delta\in\hat{K}$, define the {\it generalized spherical function} $\Phi_{s,\delta}$ on $A_1$
by
\begin{equation*}
\Phi_{s,\delta}(a_t)=(\mathcal{P}_s\varphi_\delta)(a_t).
\end{equation*}

\begin{proposition}\label{prop_Poisson_k_type}
Let $s\in\mathbb{C}$, $\delta\in\hat{K}$ and  $f\in V_\delta$. Then for any
$k\in K$ and $a_t\in A_1$,
\begin{equation*}
(\mathcal{P}_sf)(ka_t)=\Phi_{s,\delta}(a_t)f(k).
\end{equation*}
\end{proposition}
\begin{proof}
Since $M_1$ centralize $A_1$, we can view the operator (\ref{op_P_t}) as a
bounded operator on  $L^2(S)$. Moreover, $P^t_s$ commutes with the left
regular representation of $K$ in $L^2(S)$. Hence, by Schur's lemma,
$P^t_s=\Phi_{s,\delta}(a_t)\cdot I$ on each $V_\delta$ and the proposition follows.
\end{proof}

The first main theorem of this section can now be stated as follows :

\begin{theorem}\label{main_th_L_2}
Let $s\in\mathbb{C}$ be such that $\Re(s)>\frac{a}{2}(r-1)$. A smooth
function $F$ on $\Omega$ is the Poisson transform $F=\mathcal{P}_sf$ of a
function $f\in L^2(S)$ if and only if
$F\in\mathcal{E}_s^2(\Omega)$.
\end{theorem}

\begin{proof}
The necessary condition follows from Proportion \ref{nece_cond} and
\cite[Theorem 6.1 and Theorem 7.2]{Koufany-Zhang}. On the other hand,
let $F\in
\mathcal{E}_s^p(\Omega)$. We apply again \cite[Theorem 6.1 and Theorem
7.2]{Koufany-Zhang}. Then, there exists a hyperfunction
$f\in\mathcal{B}(S)$ such that $F=\mathcal{P}_sf$. Let
$f=\sum_{\delta\in\hat{K}}f_\delta$ be its $K-$type decomposition. By
Proposition \ref{prop_Poisson_k_type} we can write 
$$F(ka_t)=\sum_{\delta\in\hat{K}}\Phi_{s,\delta}(a_t)f_\delta(k)$$
in $\mathcal{C}^\infty(K\times [0,+\infty[)$. \\
Now observe that
$$\|F\|_{s,2}^2=\sup_{t>0}e^{2(n-r\Re(s))t}\sum_{\delta\in\hat{K}}|\Phi_{s,\delta}(a_t)|^2
\|f_\delta\|_2^2<\infty. $$
Then, if $\Lambda$ is an arbitrary finite subset of $\hat{K}$, we get
$$e^{2(n-r\Re(s))t}\sum_{\delta\in\Lambda}|\Phi_{s,\delta}(a_t)|^2
\|f_\delta\|_2^2\leq \|F\|_{s,2}^2$$
for every $t>0$ and hence form Theorem \ref{thm_Fatou1} it follows immediately that
$$|\mathbf{c}_s|^2 \sum_{\delta\in\Lambda}\|f_\delta\|_2^2\leq \|F\|^2_{s,2}$$
which implies that $f=\sum_{\delta\in K}f_\delta$ in $L^2(S)$ and that
$$|\mathbf{c}_s|^2\|f\|_2\leq \|F\|_{s,2}.$$
This ends the proof of the theorem.
\end{proof}

In the following proposition we show how to recover a function $f\in L^2(S)$ from its
Poisson transform $\mathcal{P}_sf$.

\begin{proposition}\label{inversion_L_2}
Let $s\in\mathbb{C}$ be such that $\Re(s)>\frac{a}{2}(r-1)$. Let
$F\in\mathcal{E}_s^2(\Omega)$ and $f\in L^2(S)$ its boundary value. Then
the following inversion formula
\begin{equation}\label{inversion_formula}
f(k)=|\mathbf{c}_s|^{-2}\lim_{t\to\infty}e^{2(n-r\Re(s))t}
\int_K\overline{e^{-\langle s\rho_0+\rho_1, H_1(a_{-t}k^{-1}h)\rangle }}F(ha_t)dh
\end{equation}
holds in $L^2(S)$.
\end{proposition}

\begin{proof}
Let $F\in\mathcal{E}_s^2(\Omega)$, then its follows from Theorem \ref{main_th_L_2} that
there exists a unique $f\in L^2(S)$ such that $F=\mathcal{P}_sf$. Let
$f=\sum_{\delta\in\hat{K}}f_\delta$ be its $K-$type expansion, then
similarly to the preceding proof, we get
\begin{equation}
F(ka_t)=\sum_{\delta\in\hat{K}}\Phi_{s,\delta}(a_t)f_\delta(k).
\end{equation}
For any $t>0$, define the complex-valued function on $K$ by
$$g_t(ka_t)=|\mathbf{c}_s|^{-2} e^{2(n-rs)t}\int_K\overline{e^{-\langle
    s\rho_0+\rho_1, H_1(a_{-t}k^{-1}h)\rangle }}F(ha_t)dh.$$
Next, using the above series expansion  we can write according to Theorem \ref{thm4.5},
$$g_t(h)=|\mathbf{c}_s|^{-2} e^{2(n-rs)t}\sum_{\delta\in\hat{K}}
|\Phi_{s,\delta}(a_t)|^2 f_\delta(h).$$
Thus,
$$
\|g_t-f\|_2^2=\sum_{\delta\in\hat{K}}\left| |\mathbf{c}_s|^{-2} e^{2(n-rs)t}
  |\Phi_{s,\delta}(a_t)|^2-1\right|^2 \|f_\delta\|_2^2, $$ 
which shows that $\|g_t-f\|_2\to 0$,
since $\lim_{t\to\infty} e^{(n-rs)t}\Phi_{s,\delta}(a_t)=\mathbf{c}_s$.
\end{proof}

\subsection{The $\mathrm{L}^p-$Poisson transform range, $p\not=2$}
We shall now prove the second main result of this paper, more precisely,
we shall characterize the $L^p-$range of the Poisson transform.
We will need the following notation. For each function $f$ on $\Omega$, define the
function $f^t$, $t>0$, on $K$ by
\begin{equation*}
f^t(k)=f(ka_t).
\end{equation*}

\begin{theorem}\label{second_main_Thm}
Let $s\in\mathbb{C}$ be such that $\Re(s)>\frac{a}{2}(r-1)$. A smooth
function $F$ on $\Omega$ is the Poisson transform $F=\mathcal{P}_sf$ of a
function $f\in L^p(S)$ if and only if $F\in\mathcal{E}_s^p(\Omega)$.
\end{theorem}

\begin{proof}
We will follow the technique used by Kor{\`a}nyi \cite{Koranyi}.  
Let $(\chi_n)_n$ be an approximation of the identity in
$\mathcal{C}(K)$. That is $\chi_n\geq 0$, $\int_K\chi_n(k)dk=1$ and
$lim_{n\to+\infty}\int_{K\setminus U}\chi_n(k)dk=0$ for every neighborhood
$U$ of $e$ in $K$. Let $F\in
\mathcal{E}_s^p(\Omega)$. For each $n$, define the function $F_n$ on $\Omega$ by
$$F_n(gK)=\int_K\chi_n(k)F(k^{-1}g)dk.$$ Then
$(F_n)_n$ converges point-wise to $F$, and since the set
$\mathcal{E}_s(\Omega)$ of harmonic functions
satisfying the Hua system is $G-$invariant,
$F_n\in\mathcal{E}_s(\Omega)$, for each $n$. Furthermore,
$$F_n^t(ka_tK)=(\chi_n \ast F^t)(k)$$
and this shows
\begin{equation}\label{inega_1}
\|F^t_n\|_2\leq \|\chi_n\|_2 \|F^t\|_p,
\end{equation}
and
\begin{equation}\label{inega_2}
\|F^t_n\|_p\leq \|F^t\|_p.
\end{equation}
It follows from (\ref{inega_1}) 
\begin{equation*}
\sup_{t>0} e^{(n-rs)t}\left(  \int_K  |F_n(ka_t)|^2 dk\right)^{1/2}\leq \|\chi_n\|_2 \|F\|_{s,p}.
\end{equation*}
Thus $F_n\in\mathcal{E}_{s,2}(\Omega)$ and by Theorem \ref{main_th_L_2}, there
exists $f_n\in L^2(S)$ such that $F_n=\mathcal{P}_s f_n$. Now, our goal is
to prove that $f_n$ belongs to $L^p(S)$. Using the inversion formula
(\ref{inversion_formula}) we can write in $L^2(S)$,
$$f_n(k)=\lim_{t\to+\infty} g_n^t(k)$$
where 
$$g_n^t(h)=g_n(ha_t)=
|\mathbf{c}_s|^{-2} e^{2(n-r\Re(s))t} \int_K \overline{e^{-\langle s\rho_0+\rho_1,H_1(a_{-t}k^{-1}h) \rangle }} F_n(ka_t)dk.$$
Let $\varphi\in\mathcal{C}(S)$ be a continuous function on $S$, then
$$\int_K f_n(h)\varphi(h)dh=\lim_{t\to \infty}\int_K g_n^t(h)\varphi(h) dh.$$
Moreover,

\begin{eqnarray*}
\int_K g_n^t(h)\varphi(h) dh
 &=& |\mathbf{c}_s|^{-2}  e^{2(n-r\Re(s))t} \times\\ 
&\times & \int_K \int_K F_n(ka_t) \varphi(h) 
 \overline{e^{-\langle s\rho_0+\rho_1,H_1(a_{-t}k^{-1}h) \rangle }} dkdh\\
 &=& |\mathbf{c}_s|^{-2}  e^{2(n-r\Re(s))t} \int_K\overline{\mathcal{P}_s\bar{\varphi}(ka_t)} F_n(ka_t)dk.
\end{eqnarray*}
By the Holder inequality, if $q$ is such that $1/p + 1/q=1$, we get
\begin{eqnarray*}
\left| \int_K g_n^t(h)\varphi(h) dh \right|
&\leq&  |\mathbf{c}_s|^{-2}  e^{2(n-r\Re(s))t} \| \mathcal{P}_s\varphi \|_q \| F_n^t\|_p,\\
&\leq&  |\mathbf{c}_s|^{-2}  e^{2(n-r\Re(s))t} \| \mathcal{P}_s\varphi \|_q \| F^t\|_p,
\end{eqnarray*}
where the second inequality follows from (\ref{inega_2}).  But $F\in\mathcal{E}_{s,p}(\omega)$, then
$$\left| \int_K g_n^t(h)\varphi(h) dh \right| \leq 
|\mathbf{c}_s|^{-2}  e^{2(n-r\Re(s))t} \| \mathcal{P}_s\varphi \|_q \| F\|_{s,p}.$$
Therefore, by Theorem \ref{thm_Fatou1},
\begin{equation*}\label{eq:f_n}
\left| \int_K f_n(h)\varphi(h)dh \right| \leq 
|\mathbf{c}_s|^{-1}  \| \varphi \|_q \| F \|_{s,p},
\end{equation*}
and by taking the supremum over $\varphi\in\mathcal{C}(S)$ with
$\|\varphi\|_q=1$, we get
$$\| f_n\|_p \leq |\mathbf{c}_s|^{-1} \|F\|_{s,p}.$$
Now, for each $\varphi\in L^q(S)$, define the functional 
\begin{equation*}
T_n(\varphi)=\int f_n(h)\varphi(h)dk.
\end{equation*}
Then it is obvious by (\ref{eq:f_n}) that 
$$| T_n)\varphi) |\leq |\mathbf{c}_s|^{-1} \| \varphi \|_q \|F\|_{s,p}$$ hence, $T_n$ is uniformly bounded operator
in $L^q(S)$ with $\sup_n \| T_n \| \leq |\mathbf{c}_s|^{-1} \|F\|_{s,p}$. Thanks to
Banach-Alaouglu-Bourbaki's theorem, there exists a subsequence of bounded
operators $(T_{n_j})_j$ which converges as $n_j\to +\infty$ to a bounded
operator $T$ in $L^q(S)$, under the $\ast-$weak topology, with $\| T \|
\leq |\mathbf{c}_s|^{-1} \|F\|_{s,p}$. Then, by the Riesz representation theorem,
there exists a unique function $f\in
L^p(S)$ such that 
$$T(\varphi)=\int_K f(h)\varphi(h)dh, \;\; \forall \varphi \in L^q(S)$$ with 
\begin{equation}\label{eq:f}
\| f\|_p\leq \| T_n \| \leq |\mathbf{c}_s|^{-1} \|F\|_{s,p}.
\end{equation}
Now, observe that $$F_{n_j}(g)=T_{n_j}( e^{-\langle s\rho_0+\rho_1,
  H_1(g^{-1}k\rangle}),$$ thus, by taking the limit as  $n\to +\infty$ we
get
$$F(g)=T(e^{-\langle s\rho_0+\rho_1,
  H_1(g^{-1}k\rangle})=\mathcal{P}_s f(g)$$
with $|\mathbf{c}_s| \|f\|_p\leq \|F\|_{s,p}$, by (\ref{eq:f}), and this finishes the
proof of the theorem.
\end{proof}
\bibliographystyle{amsplain}

\end{document}